\documentclass[11pt]{article}
\usepackage{amscd}
\usepackage{amsfonts}
\usepackage{amsmath}
\usepackage{amssymb}
\usepackage{amsthm}
\usepackage{bbm}
\usepackage{CJK}
\usepackage{fancyhdr}
\usepackage{graphicx}
\usepackage{hyperref}
\usepackage{indentfirst}
\usepackage{latexsym}
\usepackage{mathrsfs}
\usepackage{xypic}
\allowdisplaybreaks

\usepackage[top=1in,bottom=1in,left=1.25in,right=1.25in]{geometry}
\def\<{\langle}
\def\>{\rangle}
\def\a{\alpha}

\def\c{\cdot}

\def\D{\Delta}

\def\lr{\longrightarrow}

\def\o{\otimes}

\def\v{\epsilon }

\date{}
\begin{document}
\renewcommand{\baselinestretch}{1.2}
\renewcommand{\arraystretch}{1.0}
\title{\bf Hom-Yang-Baxter equations and Hom-Yang-Baxter systems}
\author{{\bf Shengxiang Wang$^{1}$, Xiaohui Zhang$^{2}$,
        Shuangjian Guo$^{3}$\footnote
        {Correspondence: shuangjianguo@126.com} }\\
1.~ School of Mathematics and Finance, Chuzhou University,\\
 Chuzhou 239000,  China \\
2.~  School of Mathematical Sciences, Qufu Normal University, \\Qufu 273165, China\\
 3.~ School of Mathematics and Statistics, Guizhou University of\\ Finance and Economics, Guiyang 550025, China}
 \maketitle
\begin{center}
\begin{minipage}{13.cm}

{\bf \begin{center} ABSTRACT \end{center}}
In this paper, we mainly present some new solutions of the Hom-Yang-Baxter equation from  Hom-algebras,
 Hom-coalgebras and  Hom-Lie  algebras, respectively.
 Also, we prove that these solutions are all self-inverse and give some examples.
  Finally, we introduce the notion of  Hom-Yang-Baxter systems and obtain two kinds of Hom-Yang-Baxter systems.\\


 \end{minipage}
 \end{center}
 \normalsize\vskip1cm

\section*{INTRODUCTION}
\def\theequation{0. \arabic{equation}}
\setcounter{equation} {0}

The Yang-Baxter equation (YBE), which was first introduced by Yang, Baxter and McGuire \cite{Baxter, yang}.
YBE has various forms   in physics and plays an important role in many topics in mathematical physics, including quantum groups, quantum integrable systems, braided categories  and invariants of knots and links.

 Recently, a twisted Hom-type generalization of the YBE called Hom-Yang-Baxter equation (HYBE) was introduced in \cite{Yau2009, Yau2011, Yau2012} by Yau.
 The HYBE states
\begin{eqnarray*}
(\a\o B)\circ(B\o \a)\circ(\a\o B)=(B\o \a)\circ(\a\o B)\circ(B\o \a),
\end{eqnarray*}
where $\a$ is an endomorphism of the vector space $V$, and $B: V\o V \rightarrow V\o V$ is a bilinear map that commutes with $\a\o \a$.
In which Yau constructed several classes of solutions of the HYBE, generalizing the solutions of the YBE from Lie algebras and  quasitriangular bialgebras.
Later, Yau  \cite{Yau2015} extended the  the classical Yang-Baxter equation to the classical Hom-Yang-Baxter equation (CHYBE) in a Hom-Lie algebra
and studied the related algebraic structure.

The study of Hom-algebras can be traced back to  Hartwig,  Larsson and  Silvestrov's work in \cite{Hartwig},
 where the notion of Hom-Lie algebra in the context of q-deformation theory of Witt and Virasoro algebras \cite{Hu} was introduced,
which plays an important role in physics, mainly in conformal field theory.
Hom-algebras and Hom-coalgebras were introduced by Makhlouf and Silvestrov  \cite{Makhlouf2008} as
generalizations of ordinary algebras and coalgebras in the following sense:
the associativity of the multiplication is replaced by
the Hom-associativity and similar for Hom-coassociativity.
They also defined the structures of Hom-bialgebras and Hom-Hopf algebras,
and described some of their properties extending properties of ordinary bialgebras and Hopf algebras
in \cite{Makhlouf2009, Makhlouf2010}.
Many more properties and structures of Hom-Hopf algebras have been developed,
see  \cite{CX14, Gohr2010, ma2014} and references cited therein.

In \cite{Yau2009, Yau2012} Yau proposed the definition of quasitriangular Hom-Hopf algebras
and showed that each quasitriangular Hom-Hopf algebra yields a solution of the Hom-Yang-Baxter equation.
Meanwhile, several classes of solutions of the  Hom-Yang-Baxter equation were constructed from different respects,
including those associated
to Hom-Lie algebras \cite{Fang, wang&GUO2020, Yau2009, Yau2011}, Drinfeld (co)doubles \cite{CWZ2014, ZGW, ZWZ} and Hom-Yetter-Drinfeld modules \cite{chen2014, LIU2014, ma2017, ma2019, Makhlouf2014, WCZ2012, YW}.

In \cite{Nichita1999, Nichita2006, Nichita2008}, Nichita,  Parashar and Popovici presented some solutions of the
YBE from associative algebras, coassociative coalgebras and  Lie  algebras respectively.
Depending on the above working, the motivation of   constructing new solutions of HYBE is natural.
The purpose of the present paper is to investigate how to construct solutions of HYBE from  Hom-algebras,
 Hom-coalgebras and  Hom-Lie  algebras, and show two kinds of Hom-Yang-Baxter systems from  Hom-algebras and Hom-coalgebras.

This paper is organized as follows.
In Section 1, we recall some basic definitions about Hom-algebras, Hom-coalgebras and  Hom-Lie  algebras.
In Section 2, we show two solutions of  HYBE from  Hom-algebras  and prove that they are self-inverse  (see Theorem 2.1 and Theorem 2.4).
In Section 3, we present two solutions of  HYBE from  Hom-coalgebras  and prove that they are self-inverse  (see Theorem 3.1 and Theorem 3.4).
In Section 4,  we obtain a solution  of  HYBE from  Hom-Lie algebras  and prove that it is self-inverse  (see Theorem 4.1).
In Section 5, we introduce the notion of  Hom-Yang-Baxter systems and present two kinds of Hom-Yang-Baxter systems (see Theorem 5.2 and Theorem 5.3).

\section{PRELIMINARIES}
\def\theequation{\arabic{section}.\arabic{equation}}
\setcounter{equation} {0}

Throughout this paper, $k$ is a fixed field.
 Unless otherwise stated, all vector spaces, algebras, modules, maps and unadorned tensor products are over $k$.
 For a coalgebra $C$, the coproduct  will be denoted by $\Delta$.
 We adopt a Sweedler's notation $\triangle(c)=c_{1}\otimes c_{2}$, for any $c\in C$, where the summation is understood.
We refer to \cite{Sweedler} for the Hopf algebra theory and terminology.

 We now recall some useful definitions in \cite{Li2014, Makhlouf2008, Makhlouf2009,  Makhlouf2010}.
\smallskip

 {\bf Definition 1.1.} A Hom-algebra is a quadruple $(A,\mu,1_A,\a)$ (abbr. $(A,\a)$), where $A$ is a $k$-linear space,
  $\mu: A\o A \lr A$ is a $k$-linear map, $1_A \in A$ and $\a$ is a endomorphism of $A$, such that
 \begin{eqnarray*}
 &(HA1)& \a(aa')=\a(a)\a(a');~~\a(1_A)=1_A,\\
 &(HA2)& \a(a)(a'a'')=(aa')\a(a'');~~a1_A=1_Aa=\a(a)
 \end{eqnarray*}
 are satisfied for $a, a', a''\in A$. Here we use the notation $\mu(a\o a')=aa'$.
\smallskip

 {\bf Definition 1.2.}  A Hom-coalgebra is a quadruple $(C,\D,\varepsilon,\a)$ (abbr. $(C,\a)$),
  where $C$ is a $k$-linear space, $\D: C \lr C\o C$, $\varepsilon: C\lr k$ are $k$-linear maps,
  and $\a$ is a endomorphism of $C$, such that
 \begin{eqnarray*}
 &(HC1)& \a(c)_1\o \a(c)_2=\a(c_1)\o \a(c_2);~~\varepsilon\circ \a=\varepsilon;\\
 &(HC2)& \a(c_{1})\o c_{21}\o c_{22}=c_{11}\o c_{12}\o \a(c_{2});~~\v(c_1)c_2=c_1\v(c_2)=\a(c)
 \end{eqnarray*}
 are satisfied for $c \in C$.
\smallskip

 {\bf Definition 1.3.}
  A Hom-Lie algebra  is a triple $(L, [\c,\c],\alpha)$ consisting of a linear space $L$,
a bilinear map $[\c,\c]:L\o L\rightarrow L$
and a endomorphism $\alpha:L\rightarrow L$, such that
\begin{eqnarray*}
&(HL1)&[l,l']=-[l',l],\\
&(HL2)&[\alpha(l),[l',l'']]+[\alpha(l'),[l'',l]]+[\alpha(l''),[l,l']]=0
\end{eqnarray*}
are satisfied for all $l,l',l''\in L$.

\section{Solutions of the HYBE from Hom-algebras}
\def\theequation{\arabic{section}.\arabic{equation}}
\setcounter{equation} {0}

In this section, we will give two kinds of solutions of the HYBE from Hom-algebras and prove that these solutions are both self-inverse.
\medskip

\noindent{\bf Theorem 2.1.}
Let  $(A,\mu, 1_A, \alpha)$  be a Hom-algebra and $\lambda,\nu\in k$.
Then
\begin{eqnarray*}
B: A\o A\rightarrow A\o A,~a\o b\mapsto \lambda ab\o 1_A+\nu 1_A\o ab-\lambda \a(a)\o \a(b)
\end{eqnarray*}
is a solution for HYBE.
\medskip

{\bf Proof.} We first show that $B$ is compatible with the twist map $\a$.
For this, we take  any $a,b\in A$ and calculate
\begin{eqnarray*}
(\a\o \a)\circ B(a\o b)
&=&(\a\o \a)(\lambda ab\o 1_A+\nu 1_A\o ab-\lambda \a(a)\o \a(b))\\
&=&\lambda\a(ab)\o 1_A+\nu 1_A\o \a(ab)-\lambda \a^{2}(a)\o \a^{2}(b),\\
B\circ(\a\o \a)(a\o b)
&=&\lambda\a(a)\a(b)\o 1_A+\nu 1_A\o \a(a)\a(b)-\lambda \a^{2}(a)\o \a^{2}(b).
\end{eqnarray*}
It follows that $(\a\o \a)\circ B=B\circ(\a\o \a)$, as desired.

 Next we will verify that $B$ satisfies the HYBE. In fact, for any  $a,b,c\in A$,
 one may directly check that
$
 B(a\o 1_A)=\nu 1_A\o\a(a).
$
 On the one hand, we have
\begin{small}
\begin{eqnarray*}
&&(\a\o B)(a\o b\o c)\\
&&~~~~=\a(a)\o (\lambda bc\o 1_A+\nu 1_A\o bc-\lambda \a(b)\o \a(c)),\\
&&(B\o \a)\circ(\a\o B)(a\o b\o c)\\
&&~~~~=\lambda B(\a(a)\o  bc)\o 1_A+\nu B(\a(a)\o 1_A)\o \a(bc)-\lambda B(\a(a)\o \a(b))\o \a^{2}(c)\\
&&~~~~=\lambda ^{2}\a(a)(bc)\o 1_A\o 1_A+\lambda\nu 1_A\o\a(a)(bc)\o 1_A-\lambda ^{2}\a^{2}(a)\o\a(bc)\o 1_A\\
&&~~~~~~~+\nu^{2}1_A\o  \a^{2}(a)\o\a(bc)-\lambda ^{2}\a(ab)\o 1_A\o\a^{2}(c)-\lambda\nu1_A\o\a(ab) \o\a^{2}(c)\\
&&~~~~~~~+\lambda ^{2}\a^{2}(a)\o\a^{2}(b)\o\a^{2}(c),\\
&&(\a\o B)\circ(B\o \a)\circ(\a\o B)(a\o b\o c)\\
&&~~~~=\lambda ^{2}\a^{2}(a)(bc)\o B(1_A\o 1_A)+\lambda\nu 1_A\o B(\a(a)(bc)\o 1_A)-\lambda ^{2}\a^{3}(a)\o B(\a(bc)\o 1_A)\\
&&~~~~~~~+\nu^{2}1_A\o  B(\a^{2}(a)\o\a(bc))-\lambda ^{2}\a^{2}(ab)\o B(1_A\o\a^{2}(c))-\lambda\nu1_A\o B(\a(ab) \o\a^{2}(c))\\
&&~~~~~~~+\lambda ^{2}\a^{3}(a)\o B(\a^{2}(b)\o\a^{2}(c))\\
&&~~~~=\underline{\lambda^{2}\nu\a^{2}(a)\a(bc)\o 1_A\o 1_A}_{(1)}+\underline{\lambda\nu^{2}1_A\o 1_A\o\a^{2}(a)\a(bc)}_{(2)}
        -\underline{\lambda^{2}\nu\a^{3}(a)\o 1_A\o\a^{2}(bc)}_{(3)}\\
&&~~~~~~~+\underline{\lambda\nu^{2}1_A\o\a^{2}(a)\a(bc)\o 1_A}_{(4)}
         +\underline{\nu^{3}1_A\o 1_A\o\a^{2}(a)\a(bc)}_{(5)}-\underline{\lambda\nu^{2}1_A\o\a^{3}(a)\o\a^{2}(bc)}_{(6)}\\
&&~~~~~~~-\underline{\lambda^{3}\a^{2}(ab)\o\a^{3}(c)\o 1_A}_{(7)}-\underline{\lambda^{2}\nu\a^{2}(ab)\o 1_A\o\a^{3}(c)}_{(8)}
         +\underline{\lambda^{3}\a^{2}(ab)\o 1_A\o\a^{3}(c)}_{(9)}\\
&&~~~~~~~-\underline{\lambda^{2}\nu 1_A\o\a(ab)\a^{23}(c)\o 1_A}_{(10)}
         -\underline{\lambda\nu^{2}1_A\o 1_A\o\a^{2}(a)\a(bc)}_{(11)}+\underline{\lambda^{2}\nu 1_A\o\a^{2}(ab)\o\a^{3}(c)}_{(12)}\\
&&~~~~~~~+\underline{\lambda^{3}\a^{3}(a)\o\a^{2}(bc)\o 1_A}_{(13)}+\underline{\lambda^{2}\nu\a^{3}(a)\o 1_A\o\a^{2}(bc)}_{(14)}
         -\underline{\lambda^{3}\a^{3}(a)\o\a^{3}(b)\o\a^{3}(c)}_{(15)}\\
&&~~~~=\underline{\lambda^{2}\nu\a^{2}(a)\a(bc)\o 1_A\o 1_A}_{(1)}+\underline{\nu^{3}1_A\o 1_A\o\a^{2}(a)\a(bc)}_{(5)}
         -\underline{\lambda\nu^{2}1_A\o\a^{3}(a)\o\a^{2}(bc)}_{(6)}\\
&&~~~~~~~-\underline{\lambda^{3}\a^{2}(ab)\o\a^{3}(c)\o 1_A}_{(7)}
         -\underline{\lambda^{2}\nu\a^{2}(ab)\o 1_A\o\a^{3}(c)}_{(8)}+\underline{\lambda^{3}\a^{2}(ab)\o 1_A\o\a^{3}(c)}_{(9)}\\
&&~~~~~~~+\underline{\lambda^{2}\nu 1_A\o\a^{2}(ab)\o\a^{3}(c)}_{(12)}+\underline{\lambda^{3}\a^{3}(a)\o\a^{2}(bc)\o 1_A}_{(13)}
         -\underline{\lambda^{3}\a^{3}(a)\o\a^{3}(b)\o\a^{3}(c)}_{(15)}.
\end{eqnarray*}
\end{small}
The last equality holds since (2)=(11), (3)=(14) and  (4)=(10).

 On the other hand, we have
\begin{small}
\begin{eqnarray*}
&&(B\o \a)(a\o b\o c)\\
&&~~~~=\lambda ab\o 1_A\o\a(c) +\nu 1_A\o ab\o\a(c)-\lambda \a(a)\o \a(b)\o \a(c),\\
&&(\a\o B)\circ(B\o \a)(a\o b\o c)\\
&&~~~~=\lambda\a(ab)\o B(1_A\o\a(c))+\nu 1_A\o B(ab\o\a(c))-\lambda \a^{2}(a)\o B(\a(b)\o \a(c))\\
&&~~~~=\lambda^{2}\a(ab)\o\a^{2}(c)\o 1_A+\lambda\nu\a(ab)\o 1_A\o\a^{2}(c)-\lambda^{2}\a(ab)\o 1_A\o\a^{2}(c)\\
&&~~~~~~~+\lambda\nu 1_A\o(ab)\a(c)\o1_A+\nu^{2}1_A\o1_A\o(ab)\a(c)-\lambda\nu 1_A\o\a(ab)\o\a^{2}(c)\\
&&~~~~~~~-\lambda^{2}\a^{2}(a)\o\a(bc)\o 1_A-\lambda\nu\a^{2}(a)\o 1_A\o\a(bc)+\lambda^{2}\a^{2}(a)\o\a^{2}(b)\o\a^{2}(c),\\
&&(B\o \a)\circ(\a\o B)\circ(B\o \a)(a\o b\o c)\\
&&~~~~=\lambda^{2}B(\a(ab)\o\a^{2}(c))\o 1_A+\lambda\nu B(\a(ab)\o 1_A)\o\a^{3}(c)-\lambda^{2} B(\a(ab)\o 1_A)\o\a^{3}(c)\\
&&~~~~~~~+\lambda\nu B(1_A\o(ab)\a(c))\o1_A+\nu^{2}B(1_A\o1_A)\o\a(ab)\a^{2}(c)-\lambda\nu B(1_A\o\a(ab))\o\a^{3}(c)\\
&&~~~~~~~-\lambda^{2}B(\a^{2}(a)\o\a(bc))\o 1_A-\lambda\nu B(\a^{2}(a)\o 1_A)\o\a^{2}(bc)+\lambda^{2}B(\a^{2}(a)\o\a^{2}(b))\o\a^{3}(c)\\
&&~~~~=\underline{\lambda^{3}\a(ab)\a^{2}(c)\o1_A\o1_A}_{(16)}+\underline{\lambda^{2}\nu 1_A\o\a(ab)\a^{2}(c)\o1_A}_{(17)}
       -\underline{\lambda^{3}\a^{2}(ab)\o\a^{3}(c)\o1_A}_{(18)}\\
&&~~~~~~~-\underline{\lambda^{2}\nu1_A\o\a^{2}(ab)\o\a^{3}(c)}_{(19)}+\underline{\lambda\nu^{2}1_A\o\a^{2}(ab)\o\a^{3}(c)}_{(20)}
         +\underline{\lambda^{2}\nu\a(ab)\a^{2}(c)\o1_A\o1_A}_{(21)}\\
&&~~~~~~~+\underline{\lambda\nu^{2}1_A\o\a(ab)\a^{2}(c)\o1_A}_{(22)}-\underline{\lambda^{2}\nu 1_A\o\a(ab)\a^{2}(c)\o1_A}_{(23)}
         +\underline{\nu^{3}1_A\o1_A\o\a(ab)\a^{2}(c)}_{(24)}\\
&&~~~~~~~-\underline{\lambda^{2}\nu\a^{2}(ab)\o 1_A\o\a^{3}(c)}_{(25)}-\underline{\lambda\nu^{2}1_A\o\a^{2}(ab)\o\a^{3}(c)}_{(26)}
         +\underline{\lambda^{2}\nu1_A\o\a^{2}(ab)\o\a^{3}(c)}_{(27)}\\
&&~~~~~~~-\underline{\lambda^{3}\a^{2}(a)\a(bc)\o1_A\o1_A}_{(28)}-\underline{\lambda\nu^{2}1_A\o\a^{2}(a)\a(bc)\o1_A}_{(29)}
         +\underline{\lambda^{3}\a^{3}(a)\o\a^{2}(bc)\o1_A}_{(30)}\\
&&~~~~~~~-\underline{\lambda\nu^{2}1_A\o\a^{3}(a)\o\a^{2}(bc)}_{(31)}+\underline{\lambda^{3}\a^{2}(ab)\o1_A\o\a^{3}(c)}_{32}
         +\underline{\lambda^{2}\nu 1_A\o\a^{2}(ab)\o\a^{3}(c)}_{(33)}\\
&&~~~~~~~-\underline{\lambda^{3}\a^{3}(a)\o\a^{3}(b)\o\a^{3}(c)}_{(34)}\\
&&~~~~=-\underline{\lambda^{3}\a^{2}(ab)\o\a^{3}(c)\o1_A}_{(18)}+\underline{\lambda^{2}\nu\a(ab)\a^{2}(c)\o1_A\o1_A}_{(21)}
        +\underline{\nu^{3}1_A\o1_A\o\a(ab)\a^{2}(c)}_{(24)}\\
&&~~~~~~~-\underline{\lambda^{2}\nu\a^{2}(ab)\o 1_A\o\a^{3}(c)}_{(25)}+\underline{\lambda^{2}\nu1_A\o\a^{2}(ab)\o\a^{3}(c)}_{(27)}
         +\underline{\lambda^{3}\a^{3}(a)\o\a^{2}(bc)\o1_A}_{(30)}\\
&&~~~~~~~-\underline{\lambda\nu^{2}1_A\o\a^{3}(a)\o\a^{2}(bc)}_{(31)}+\underline{\lambda^{3}\a^{2}(ab)\o1_A\o\a^{3}(c)}_{(32)}
        -\underline{\lambda^{3}\a^{3}(a)\o\a^{3}(b)\o\a^{3}(c)}_{(34)}.
\end{eqnarray*}
\end{small}
The last equality holds since (16)=(28), (17)=(23), (19)=(33), (20)=(26) and  (22)=(29).

Compare two expressions above, we have
 \begin{eqnarray*}
 &&(1)=(21), (5)=(24), (6)=(31), (7)=(18),\\
 &&(8)=(25), (9)=(32), (12)=(27), (13)=(30), (15)=(34).
 \end{eqnarray*}
It follows that
\begin{eqnarray*}
(\a\o B)\circ(B\o \a)\circ(\a\o B)(a\o b\o c)
=(B\o \a)\circ(\a\o B)\circ(B\o \a)(a\o b\o c).
\end{eqnarray*}
So $B$ is a solution for HYBE. $\hfill \Box$
\medskip

\noindent{\bf Corollary 2.2.}
Let  $(A,\mu, 1_A, \alpha)$  be a Hom-algebra and $\lambda, \nu\in k^{\ast}$.
Assume that $\a$ is involutive, then the solution $B$ in Theorem 2.1 is   invertible,
where the inverse is given by
\begin{eqnarray*}
B^{-1}: A\o A\rightarrow A\o A,~a\o b\mapsto \frac{1}{\nu }ab\o 1_A+\frac{1}{\lambda}1_A\o ab-\frac{1}{\lambda} \a(a)\o \a(b).
\end{eqnarray*}

{\bf Proof.}
We first show that $B\circ B^{-1}=id_{A\o A}$. In fact, for any  $a,b,c\in A$,
we have
\begin{eqnarray*}
B\circ B^{-1}(a\o b)
&=&\frac{1}{\nu }B(ab\o 1_A)+\frac{1}{\lambda}B(1_A\o ab)-\frac{1}B({\lambda} \a(a)\o \a(b))\\
&=&1_A\o\a(ab)+\frac{1}{\lambda}\{\lambda\a(ab)\o 1_A+\nu1_A\o\a(ab)-\lambda1_A\o\a(ab)\}\\
   &&-\frac{1}{\lambda}\{\lambda\a(ab)\o 1_A+\nu1_A\o\a(ab)-\lambda\a^{2}(a)\o\a^{2}(b)\}\\
&=&a\o b.
\end{eqnarray*}
It follows that $B\circ B^{-1}=id_{A\o A}$.
 Similarly, one may check that $B^{-1}\circ B=id_{A\o A}$.  So  $B^{-1}$ is the inverse of $B$.
$\hfill \Box$
\medskip

\noindent{\bf Example 2.3.}
Let $\{x_{1},x_{2},x_{3}\}$ be a basis of a 3-dimensional linear space $A$. The
following multiplication $\mu$ and the twist map $\alpha$ on $A$ define a Hom-algebra:
\begin{eqnarray*}
&&\mu(x_{1},x_{1})=x_{1},~\mu(x_{1},x_{2})=x_{2},~\mu(x_{1},x_{3})=lx_{3},\\
&&\mu(x_{2},x_{1})=x_{2},~\mu(x_{2},x_{2})=x_{2},~\mu(x_{2},x_{3})=lx_{3},\\
&&\mu(x_{3},x_{1})=lx_{3},~\mu(x_{3},x_{2})=0,~\mu(x_{3},x_{3})=0,\\
&&\alpha(x_{1})=x_{1},~\alpha(x_{2})=x_{2},~\alpha(x_{3})=lx_{3},
\end{eqnarray*}
where $l$ is a parameter  in $k$ (\cite{ma2017}).
It is easy to see that $1_A=x_1$.
Therefore, by Theorem 2.1,  the solution $B$ of the HYBE for the Hom-algebra $A$ satisfies
 \begin{eqnarray*}
&&B(x_{1}\o x_{1})=\nu x_{1}\o x_{1},
B(x_{1}\o x_{2})=\lambda x_{2}\o x_{1}+(\nu-\lambda) x_{1}\o x_{2},\\
&&B(x_{1}\o x_{3})=\lambda l x_{3}\o x_{1}+l(\nu-\lambda) x_{1}\o x_{3},
B(x_{2}\o x_{1})=\nu x_{1}\o x_{2},\\
&&B(x_{2}\o x_{2})=\lambda x_{2}\o x_{1}+\nu x_{1}\o x_{2}-\lambda x_{2}\o x_{2},\\
&&B(x_{2}\o x_{3})=\lambda l x_{3}\o x_{1}+ \nu l x_{1}\o x_{3}-\lambda l x_{2}\o x_{3},\\
&&B(x_{3}\o x_{1})=\nu l x_{1}\o x_{3},
B(x_{3}\o x_{2})=-\lambda l x_{3}\o x_{2},
B(x_{3}\o x_{3})=-\lambda l^{2} x_{3}\o x_{3},
  \end{eqnarray*}
where $\lambda,\nu$ are two parameters in $k$.
\medskip

\noindent{\bf Theorem 2.4.}
Let  $(A,\mu, 1_A, \alpha)$  be a Hom-algebra and $\lambda,\nu\in k$.
Then
\begin{eqnarray*}
B: A\o A\rightarrow A\o A,~a\o b\mapsto \lambda ab\o 1_A+\nu 1_A\o ab-\nu \a(a)\o \a(b)
\end{eqnarray*}
is also a solution for HYBE.
If, in addition, $\lambda, \nu\in k^{\ast}$ and $\a$ is involutive, then  $B$  is   invertible,
where the inverse is given by
\begin{eqnarray*}
B^{-1}: A\o A\rightarrow A\o A,~a\o b\mapsto \frac{1}{\nu }ab\o 1_A+\frac{1}{\lambda}1_A\o ab-\frac{1}{\nu} \a(a)\o \a(b).
\end{eqnarray*}

{\bf Proof.} Similar to the proof of Theorem 2.1 and Corollary 2.2.$\hfill \Box$
\medskip

\noindent{\bf Example 2.5.}
Let $\{1,g,x,y\}$ be a basis of a 4-dimensional linear space $H_{4}$. The
following multiplication $\mu$ and the twist map $\alpha$ on $H_{4}$ define a Hom-algebra:
\begin{eqnarray*}
&&\mu(1,1)=1,~\mu(1,g)=g,~\mu(1,x)=kx,~\mu(1,y)=ky,\\
&&\mu(g,1)=g,~\mu(g,g)=g,~\mu(g,x)=ky,~\mu(g,y)=kx,\\
&&\mu(x,1)=kx,~\mu(x,g)=ky,~\mu(x,x)=0,~\mu(x,y)=0,\\
&&\mu(y,1)=ky,~\mu(y,g)=-kx,~\mu(y,x)=0,~\mu(y,y)=0,\\
&&\alpha(1)=1,~\alpha(g)=g,~\alpha(x)=kx,~\alpha(y)=ky,
\end{eqnarray*}
where $k$ is a parameter  in $k$ (\cite{Lu}).
By Theorem 2.4,  the solution $B$ of the HYBE for the Hom-algebra $H_{4}$ satisfies
 \begin{eqnarray*}
&&B(1\o 1)=\lambda 1\o 1,
B(1\o g)=\lambda g\o 1,
B(1\o x)=\lambda k x\o 1,
B(1\o y)=\lambda k 1\o y,\\
&&B(g\o 1)=(\lambda-\nu) g\o 1+\nu 1\o g,
B(g\o g)=(\lambda+\nu)1\o 1 -\nu g\o g,\\
&&B(g\o x)=\lambda k y\o 1+\nu k 1\o y-\nu k g\o x,
B(g\o y)=\lambda k x\o 1+\nu k 1\o x-\nu k g\o y,\\
&&B(x\o 1)=k(\lambda-\nu) x\o 1+\nu k 1\o x,
B(x\o g)=-\lambda k y\o 1-\nu k 1\o y-\lambda k x\o g,\\
&&B(x\o x)=-k^{2} x\o x,
B(x\o y)=-k^{2} x\o y,
B(y\o 1)=k(\lambda-\nu) y\o 1+\nu k 1\o y,\\
&&B(y\o g)=-\lambda k x\o 1-\nu k 1\o x-\lambda k y\o g,
B(y\o x)=-k^{2} y\o x,
B(y\o y)=-k^{2} y\o y,
  \end{eqnarray*}
where $\lambda,\nu$ are two parameters in $k$.

\section{Solutions of the HYBE from Hom-coalgebras}
\def\theequation{\arabic{section}.\arabic{equation}}
\setcounter{equation} {0}

In this section,   we will show two kinds of solutions of the HYBE from Hom-coalgebras and prove that these solutions are both self-inverse.
 \medskip

\noindent{\bf Theorem 3.1.}
Let  $(C,\Delta, \varepsilon, \alpha)$  be a Hom-coalgebra and $\lambda,\nu\in k$.
Then
\begin{eqnarray*}
B: C\o C\rightarrow C\o C,~a\o b\mapsto \lambda \varepsilon(a)b_1\o b_2+\nu \varepsilon(b)a_1\o a_2-\lambda \a(a)\o \a(b)
\end{eqnarray*}
is a solution for HYBE.
\medskip

{\bf Proof.} We first show that $B$ is compatible with the twist map $\a$.
For this, we take  any $a,b\in C$ and calculate
\begin{small}
\begin{eqnarray*}
(\a\o \a)\circ B(a\o b)
&=&(\a\o \a)(\lambda \varepsilon(a)b_1\o b_2+\nu \varepsilon(b)a_1\o a_2-\lambda \a(a)\o \a(b))\\
&=&\lambda \varepsilon(a)\a(b_1)\o \a(b_2)+\nu \varepsilon(b)\a(a_1)\o \a(a_2)-\lambda \a^{2}(a)\o \a^{2}(b),\\
B\circ(\a\o \a)(a\o b)
&=&\lambda \varepsilon(\a(a))\a(b_1)\o \a(b_2)+\nu \varepsilon(\a(b))\a(a_1)\o \a(a_2)-\lambda \a^{2}(a)\o \a^2(b)\\
&=&\lambda \varepsilon(a)\a(b_1)\o \a(b_2)+\nu \varepsilon(b)\a(a_1)\o \a(a_2)-\lambda \a^{2}(a)\o \a^{2}(b).
\end{eqnarray*}
\end{small}
It follows that $(\a\o \a)\circ B=B\circ(\a\o \a)$, as desired.

 Next we will verify that $B$ satisfies the HYBE. For this, we take any  $a,b,c\in C$ and calculate
\begin{small}
 \begin{eqnarray*}
&&(\a\o B)(a\o b\o c)\\
&&~~~~=\lambda\varepsilon(b) \a(a)\o c_1\o c_2+\nu \varepsilon(c)\a(a)\o b_1\o b_2-\lambda \a(a)\o \a(b)\o \a(c)),\\
&&(B\o \a)\circ(\a\o B)(a\o b\o c)\\
&&~~~~=\lambda\varepsilon(b) B(\a(a)\o c_1)\o \a(c_2)+\nu \varepsilon(c)B(\a(a)\o b_1)\o \a(b_2)-\lambda B(\a(a)\o \a(b))\o \a^{2}(c)\\
&&~~~~=\lambda^{2}\varepsilon(a)\varepsilon(b) c_{11}\o c_{12}\o \a(c_2)+\lambda\nu\varepsilon(b)\a(a_1)\o \a(a_2)\o \a^{2}(c)\\
&&~~~~~~~-\lambda^{2}\varepsilon(b) \a^{2}(a)\o \a(c_1)\o \a(c_2)+\lambda\nu\varepsilon(a)\varepsilon(c)b_{11}\o b_{12}\o \a(b_2)\\
&&~~~~~~~+\nu^{2}\varepsilon(c)\a(a_1)\o \a(a_2)\o \a^{2}(b)-\lambda\nu\varepsilon(c) \a^{2}(a)\o\a(b_1)\o \a(b_2)\\
&&~~~~~~~-\lambda^{2}\varepsilon(a)\a(b_1)\o \a(b_2)\o \a^{2}(c)-\lambda\nu\varepsilon(b)\a(a_1)\o \a(a_2)\o \a^{2}(c)\\
&&~~~~~~~+\lambda^{2}\a^{2}(a)\o \a^{2}(b)\o \a^{2}(c),\\
&&(\a\o B)\circ(B\o \a)\circ(\a\o B)(a\o b\o c)\\
&&~~~~=\lambda^{2}\varepsilon(a)\varepsilon(b) \a(c_{11})\o B(c_{12}\o \a(c_2))+\lambda\nu\varepsilon(b)\a^{2}(a_1)\o B(\a(a_2)\o \a^{2}(c))\\
&&~~~~~~~-\lambda^{2}\varepsilon(b) \a^{3}(a)\o B(\a(c_1)\o \a(c_2))+\lambda\nu\varepsilon(a)\varepsilon(c)\a(b_{11})\o B(b_{12}\o \a(b_2))\\
&&~~~~~~~+\nu^{2}\varepsilon(c)\a^{2}(a_1)\o B(\a(a_2)\o \a^{2}(b))-\lambda\nu\varepsilon(c) \a^{3}(a)\o B(\a(b_1)\o \a(b_2))\\
&&~~~~~~~-\lambda^{2}\varepsilon(a)\a^{2}(b_1)\o B(\a(b_2)\o \a^{2}(c))-\lambda\nu\varepsilon(b)\a^{2}(a_1)\o B(\a(a_2)\o \a^{2}(c))\\
&&~~~~~~~+\lambda^{2}\a^{3}(a)\o B(\a^{2}(b)\o \a^{2}(c))\\
&&~~~~=\underline{\lambda^{3}\varepsilon(a)\varepsilon(b) \a^{2}(c_{1})\o  \a(c_{21})\o  \a(c_{22})}_{(1)}
        +\underline{\lambda^{2}\nu\varepsilon(a)\varepsilon(b) \a^{2}(c_{1})\o  \a(c_{21})\o  \a(c_{22})}_{(2)}\\
&&~~~~~~~-\underline{\lambda^{2}\nu\varepsilon(a)\varepsilon(b) \a^{2}(c_{11})\o  \a(c_{12})\o  \a^{2}(c_{2})}_{(3)}
        +\underline{\lambda^{2}\nu\varepsilon(b)\a^{3}(a)\o\a^{2}(c_{1})\o  \a^{2}(c_{2})}_{(4)}\\
&&~~~~~~~+\underline{\lambda\nu^{2}\varepsilon(b)\varepsilon(c)\a^{2}(a_{1})\o  \a(a_{21})\o  \a(a_{22}) }_{(5)}
        -\underline{\lambda^{2}\nu\varepsilon(b)\a^{2}(a_1)\o\a^{2}(a_{2})\o  \a^{3}(c)}_{(6)}\\
&&~~~~~~~-\underline{\lambda^{3}\varepsilon(b) \a^{3}(a)\o\a^{2}(c_{1})\o  \a^{2}(c_{2}) }_{(7)}
        - \underline{\lambda^{2}\nu\varepsilon(b)  \a^{3}(a)\o\a^{2}(c_{1})\o  \a^{2}(c_{2})}_{(8)}  \\
&&~~~~~~~+\underline{\lambda^{3}\varepsilon(b)\a^{3}(a)\o\a^{2}(c_{1})\o  \a^{2}(c_{2})}_{(9)}
         +\underline{\lambda^{2}\nu\varepsilon(a)\varepsilon(c)\a^{2}(b_{1})\o  \a(b_{21})\o  \a(b_{22})} _{(10)}\\
&&~~~~~~~+\underline{\lambda\nu^{2}\varepsilon(a)\varepsilon(c)\a^{2}(b_{1})\o  \a(b_{21})\o  \a(b_{22}) }_{(11)}
         - \underline{\lambda^{2}\nu\varepsilon(a)\varepsilon(c)   \a(b_{11})\o  \a(b_{12})\o  \a^{2}(b_{2})}_{(12)} \\
&&~~~~~~~+\underline{\lambda\nu^{2}\varepsilon(c) \a^{3}(a)\o\a^{2}(b_{1})\o  \a^{2}(b_{2}) }_{(13)}
         +\underline{\nu^{3} \varepsilon(b)\varepsilon(c)  \a^{2}(a_1)\o\a(a_{21})\o  \a(a_{22})}_{(14)} \\
&&~~~~~~~- \underline{\lambda\nu^{2}\varepsilon(c)\a^{2}(a_{1})\o  \a^{2}(a_{2})\o \a^{3}(b) }_{(15)}
         - \underline{\lambda^{2} \nu\varepsilon(c)\a^{3}(a)\o  \a^{2}(b_{1})\o \a^{2}(b_2)}_{(16)} \\
&&~~~~~~~-\underline{\lambda\nu^{2}\varepsilon(c)\a^{3}(a)\o  \a^{2}(b_{1})\o \a^{2}(b_2)  }_{(17)}
         + \underline{ \lambda^{2} \nu\varepsilon(c)\a^{3}(a)\o  \a^{2}(b_{1})\o \a^{2}(b_2)}_{(18)}\\
&&~~~~~~~- \underline{\lambda^{3} \varepsilon(a) \a^{3}(b)\o  \a^{2}(c_{1})\o \a^{2}(c_2) }_{(19)}
         -   \underline{\lambda^{2} \nu  \varepsilon(a) \varepsilon(c)  \a^{2}(b_1)\o\a(b_{21})\o  \a(b_{22})} _{(20)}\\
&&~~~~~~~+  \underline{\lambda^{3} \varepsilon(a)  \a^{2}(b_1)\o\a^{2}(b_{2})\o\a^{3}(c)}_{(21)}
          -  \underline{\lambda^{2} \nu  \varepsilon(b) \a^{3}(a)\o  \a^{2}(c_{1})\o \a^{2}(c_2)} _{(22)}\\
&&~~~~~~~-\underline{\lambda\nu^{2}\varepsilon(b)\varepsilon(c)  \a^{2}(a_1)\o\a(a_{21})\o  \a(a_{22}) }_{(23)}
         + \underline{\lambda^{2} \nu  \varepsilon(b)\a^{2}(a_{1})\o  \a^{2}(a_{2})\o \a^{3}(c)}_{(24)}\\
&&~~~~~~~+\underline{\lambda^{3} \varepsilon(b) \a^{3}(a)\o\a^{2}(c_{1})\o \a^{2}(c_2)}_{(25)}
         + \underline{\lambda^{2} \nu \varepsilon(c) \a^{3}(a)\o\a^{2}(b_{1})\o \a^{2}(b_2)}_{(26)} \\
 &&~~~~~~~-\underline{\lambda^{3}\a^{3}(a)\o\a^{3}(b)\o\a^{3}(c)}_{(27)}\\
 &&~~~~=\underline{\lambda^{3}\varepsilon(a)\varepsilon(b) \a^{2}(c_{1})\o  \a(c_{21})\o  \a(c_{22})}_{(1)}
         -\underline{\lambda^{2}\nu\varepsilon(b)  \a^{3}(a)\o\a^{2}(c_{1})\o  \a^{2}(c_{2})}_{(8)}\\
 &&~~~~~~~+\underline{\lambda^{3} \varepsilon(b) \a^{3}(a)\o\a^{2}(c_{1})\o \a^{2}(c_2)}_{(25)}
         +\underline{\nu^{3} \varepsilon(b)\varepsilon(c)  \a^{2}(a_1)\o\a(a_{21})\o  \a(a_{22})}_{(14)}\\
 &&~~~~~~~+\underline{\lambda\nu^{2}\varepsilon(a)\varepsilon(c)\a^{2}(b_{1})\o  \a(b_{21})\o  \a(b_{22})}_{(11)}
           -\underline{\lambda^{2} \nu  \varepsilon(a) \varepsilon(c)  \a^{2}(b_1)\o\a(b_{21})\o  \a(b_{22})}_{(20)}\\
 &&~~~~~~~-\underline{\lambda\nu^{2}\varepsilon(c)\a^{2}(a_{1})\o  \a^{2}(a_{2})\o \a^{3}(b)}_{(15)}
           +\underline{\lambda^{2} \nu \varepsilon(c) \a^{3}(a)\o\a^{2}(b_{1})\o \a^{2}(b_2)}_{(26)}\\
 &&~~~~~~~-\underline{\lambda^{3} \varepsilon(a) \a^{3}(b)\o  \a^{2}(c_{1})\o \a^{2}(c_2)}_{(19)}
          +\underline{\lambda^{3} \varepsilon(a)  \a^{2}(b_1)\o\a^{2}(b_{2})\o\a^{3}(c)}_{(21)}\\
&&~~~~~~~- \underline{\lambda^{3}\a^{3}(a)\o\a^{3}(b)\o\a^{3}(c)}_{(27)}.
 \end{eqnarray*}
 \end{small}
 The last equality holds since (2)=(3), (4)=(22), (5)=(23), (6)=(24), (7)=(9), (10)=(12), (13)=(17) and  (16)=(18).

On the other hand, we have
\begin{small}
\begin{eqnarray*}
&&(B\o \a)(a\o b\o c)\\
&&~~~~=\lambda \varepsilon(a)b_1\o b_2\o \a(c)+\nu \varepsilon(b)a_1\o a_2\o \a(c)-\lambda \a(a)\o \a(b)\o \a(c),\\
&&(\a\o B)\circ(B\o \a)(a\o b\o c)\\
&&~~~~=\lambda \varepsilon(a)\a(b_1)\o B(b_2\o \a(c))+\nu \varepsilon(b)\a(a_1)\o B( a_2\o \a(c))-\lambda \a^{2}(a)\o B(\a(b)\o \a(c))\\
&&~~~~=\lambda^{2} \varepsilon(a)\a^{2}(b)\o \a(c_1)\o \a(c_2)
       +\lambda\nu\varepsilon(a)\varepsilon(c)\a(b_1)\o b_{21} \o   b_{22}\\
&&~~~~~~~-\lambda^{2} \varepsilon(a)\a(b_1)\o\a(b_2)\o\a^{2}(c)
       +\lambda\nu\varepsilon(b)\a^{2}(a)\o \a(c_1)\o \a(c_2)\\
&&~~~~~~~+\nu^{2}\varepsilon(b)\varepsilon(c)\a(a_1)\o a_{21} \o   a_{22}
       -\lambda\nu\varepsilon(b) \a(a_1)\o \a(a_2)\o\a^{2}(c)\\
&&~~~~~~~-\lambda^{2} \varepsilon(b)\a^{2}(a)\o \a(c_1)\o \a(c_2)
        -\lambda\nu\varepsilon(c)\a^{2}(a)\o \a(b_1)\o \a(b_2)\\
&&~~~~~~~+ \lambda^{2}\a^{2}(a)\o\a^{2}(b)\o\a^{2}(c),\\
&&(B\o \a)\circ (\a\o B)\circ(B\o \a)(a\o b\o c)\\
&&~~~~=\lambda^{2} \varepsilon(a)B(\a^{2}(b)\o\a(c_1))\o \a^{2}(c_2)
       +\lambda\nu\varepsilon(a)\varepsilon(c)B(\a(b_1)\o b_{21}) \o  \a(b_{22})\\
&&~~~~~~~-\lambda^{2} \varepsilon(a)B(\a(b_1)\o\a(b_2))\o\a^{3}(c)
       +\lambda\nu\varepsilon(b)B(\a^{2}(a)\o \a(c_1))\o \a^{2}(c_2)\\
&&~~~~~~~+\nu^{2}\varepsilon(b)\varepsilon(c)B(\a(a_1)\o a_{21}) \o  \a( a_{22})
       -\lambda\nu\varepsilon(b) B(\a(a_1)\o \a(a_2))\o\a^{3}(c)\\
&&~~~~~~~-\lambda^{2} \varepsilon(b)B(\a^{2}(a)\o \a(c_1))\o \a^{2}(c_2)
        -\lambda\nu\varepsilon(c)B(\a^{2}(a)\o \a(b_1))\o \a^{2}(b_2)\\
&&~~~~~~~+ \lambda^{2}B(\a^{2}(a)\o\a^{2}(b))\o\a^{3}(c)\\
&&~~~~=\underline{\lambda^{3} \varepsilon(a) \varepsilon(b)\a(c_{11})\o\a(c_{12})\o \a^{2}(c_2)} _{(28)}
       + \underline{\lambda^{2}\nu\varepsilon(a)\a^{2}(b_1) \o \a^{2}(b_2)\o \a^{3}(c)}_{(29)} \\
&&~~~~~~~-\underline{\lambda^{3} \varepsilon(a)\a^{3}(b) \o \a^{2}(c_1)\o \a^{2}(c_2) }_{(30)}
       +\underline{\lambda^{2}\nu\varepsilon(a) \varepsilon(c) \a(b_{11})\o\a(b_{12})\o \a^{2}(b_2)}_{(31)} \\
&&~~~~~~~+ \underline{\lambda\nu^{2} \varepsilon(a) \varepsilon(c) \a(b_{11})\o\a(b_{12})\o \a^{2}(b_2)} _{_{(32)}}
        -\underline{\lambda^{2}\nu\varepsilon(a) \varepsilon(c) \a^{2}(b_1) \o \a(b_{21})\o \a(b_{22})}_{(33)}\\
&&~~~~~~~-\underline{\lambda^{3} \varepsilon(a) \a^{2}(b_1) \o \a^{2}(b_2)\o \a^{3}(c)}_{(34)}
         -\underline{\lambda^{2}\nu\varepsilon(a)\a^{2}(b_1) \o \a^{2}(b_2)\o \a^{3}(c)}_{(35)}\\
&&~~~~~~~+\underline{\lambda^{3}\varepsilon(a)\a^{2}(b_1) \o \a^{2}(b_2)\o \a^{3}(c)}_{(36)}
          +\underline{\lambda^{2}\nu\varepsilon(a)\varepsilon(b)\a(c_{11})\o\a(c_{12})\o \a^{2}(c_2)}_{(37)}  \\
&&~~~~~~~+\underline{ \lambda\nu^{2} \varepsilon(b)\a^{2}(a_1) \o \a^{2}(a_2)\o \a^{3}(c) } _{(38)}
          - \underline{\lambda^{2}\nu\varepsilon(b) \a^{3}(a) \o  \a^{2}(c_1)\o \a^{2}(c_2)}_{(39)}\\
&&~~~~~~~+ \underline{\lambda\nu^{2} \varepsilon(b)\varepsilon(c)\a(a_{11})\o\a(a_{12})\o \a^{2}(a_2)}_{(40)}
         + \underline{\nu^{3} \varepsilon(b)\varepsilon(c) \a(a_{11})\o\a(a_{12})\o \a^{2}(a_2)}_{(41)}\\
&&~~~~~~~- \underline{\lambda\nu^{2} \varepsilon(b)\varepsilon(c)\a^{2}(a_1)\o \a(a_{21})\o\a(a_{22})}_{(42)}
         -\underline{\lambda^{2}\nu \varepsilon(b) \a^{2}(a_1) \o \a^{2}(a_2)\o \a^{3}(c)}_{(43)} \\
&&~~~~~~~- \underline{\lambda\nu^{2} \varepsilon(b) \a^{2}(a_1) \o \a^{2}(a_2)\o \a^{3}(c) }_{(44)}
          +\underline{\lambda^{2}\nu \varepsilon(b) \a^{2}(a_1) \o \a^{2}(a_2)\o \a^{3}(c) }_{(45)}\\
&&~~~~~~~- \underline{\lambda^{3}\varepsilon(a)\varepsilon(b)\a(c_{11})\o\a(c_{12})\o \a^{2}(c_2)}_{(46)}
         -\underline{\lambda^{2}\nu\varepsilon(b) \a^{2}(a_1) \o \a^{2}(a_2)\o \a^{3}(c)}_{(47)} \\
&&~~~~~~~+ \underline{\lambda^{3} \varepsilon(b) \a^{3}(a) \o  \a^{2}(c_1)\o \a^{2}(c_2) }_{(48)}
         -\underline{\lambda^{2}\nu \varepsilon(a)\varepsilon(c) \a(b_{11})\o\a(b_{12})\o \a^{2}(b_2)}_{(49)}\\
&&~~~~~~~-\underline{\lambda\nu^{2}\varepsilon(c) \a^{2}(a_1) \o \a^{2}(a_2)\o \a^{3}(b) } _{(50)}
         + \underline{\lambda^{2}\nu\varepsilon(c)\a^{3}(a) \o  \a^{2}(b_1)\o \a^{2}(b_2)}_{(51)}\\
&&~~~~~~~+ \underline{\lambda^{3} \varepsilon(a) \a^{2}(b_1) \o \a^{2}(b_2)\o \a^{3}(c) } _{(52)}
         + \underline{\lambda^{2}\nu \varepsilon(b) \a^{2}(a_1) \o \a^{2}(a_2)\o \a^{3}(c)}_{(53)} \\
&&~~~~~~~- \underline{\lambda^{3}  \a^{3}(a) \o\a^{3}(b) \o\a^{3}(c)}_{(54)}\\
&&~~~~=\underline{\lambda^{3} \varepsilon(a) \varepsilon(b)\a(c_{11})\o\a(c_{12})\o \a^{2}(c_2)} _{(28)}
       -\underline{\lambda^{3} \varepsilon(a)\a^{3}(b) \o \a^{2}(c_1)\o \a^{2}(c_2) }_{(30)}\\
&&~~~~~~~+ \underline{\lambda\nu^{2} \varepsilon(a) \varepsilon(c) \a(b_{11})\o\a(b_{12})\o \a^{2}(b_2)} _{_{(32)}}
        - \underline{\lambda^{2}\nu\varepsilon(b) \a^{3}(a) \o  \a^{2}(c_1)\o \a^{2}(c_2)}_{(39)}\\
&&~~~~~~~+ \underline{\nu^{3} \varepsilon(b)\varepsilon(c) \a(a_{11})\o\a(a_{12})\o \a^{2}(a_2)}_{(41)}
          + \underline{\lambda^{3} \varepsilon(b) \a^{3}(a) \o  \a^{2}(c_1)\o \a^{2}(c_2) }_{(48)}\\
&&~~~~~~~-\underline{\lambda^{2}\nu \varepsilon(a)\varepsilon(c) \a(b_{11})\o\a(b_{12})\o \a^{2}(b_2)}_{(49)}
         -\underline{\lambda\nu^{2}\varepsilon(c) \a^{2}(a_1) \o \a^{2}(a_2)\o \a^{3}(b) } _{(50)}\\
&&~~~~~~~+ \underline{\lambda^{2}\nu\varepsilon(c)\a^{3}(a) \o  \a^{2}(b_1)\o \a^{2}(b_2)}_{(51)}
          + \underline{\lambda^{3} \varepsilon(a) \a^{2}(b_1) \o \a^{2}(b_2)\o \a^{3}(c) } _{(52)}   \\
&&~~~~~~~- \underline{\lambda^{3}  \a^{3}(a) \o\a^{3}(b) \o\a^{3}(c)}_{(54)}.
\end{eqnarray*}
\end{small}
 The last equality holds since (29)=(35), (31)=(33), (34)=(36), (37)=(46), (38)=(44), (40)=(42), (43)=(45) and  (47)=(43).
Compare the two expressions above, we have
 \begin{eqnarray*}
 &&(1)=(28), (8)=(39), (11)=(32), (14)=(41), (15)=(50),\\
 &&(19)=(30), (20)=(49), (21)=(52), (25)=(48), (26)=(51), (27)=(54).
 \end{eqnarray*}
That is, we prove
 \begin{eqnarray*}
(\a\o B)\circ(B\o \a)\circ(\a\o B)(a\o b\o c)
=(B\o \a)\circ(\a\o B)\circ(B\o \a)(a\o b\o c),
\end{eqnarray*}
as desired.
So $B$ is a solution for HYBE. $\hfill \Box$
\medskip

{\bf Corollary 3.2.}
Let  $(C,\Delta, \varepsilon, \alpha)$  be a Hom-coalgebra and $\lambda,\nu\in k^{\ast}$.
Assume that $\a$ is involutive, then the solution  $B$ in Theorem 3.1 is   invertible,
where the inverse is given by
\begin{eqnarray*}
B^{-1}: C\o C\rightarrow C\o C,~a\o b\mapsto \frac{1}{\nu} \varepsilon(a)b_1\o b_2+\frac{1}{\lambda}  \varepsilon(b)a_1\o a_2-\frac{1}{\lambda} \a(a)\o \a(b)
\end{eqnarray*}

{\bf Proof.} We first show that $B\circ B^{-1}=id_{C\o C}$. In fact, for any  $a,b,c\in C$,
we have
\begin{eqnarray*}
B\circ B^{-1}(a\o b)
&=& \frac{1}{\nu} \varepsilon(a)B(b_1\o b_2)+\frac{1}{\lambda}  \varepsilon(b)B(a_1\o a_2)-\frac{1}{\lambda} B(\a(a)\o \a(b))\\
&=&\frac{1}{\nu} \varepsilon(a)\{\lambda\varepsilon(b_1)b_{21}\o b_{22}+\nu\varepsilon(b_2)b_{11}\o b_{12}-\lambda\a(b_1)\o \a(b_2)\}\\
&&+\frac{1}{\lambda}  \varepsilon(b)\{\lambda\varepsilon(a_1)a_{21}\o a_{22}+\nu\varepsilon(a_2)a_{11}\o a_{12}-\lambda\a(a_1)\o \a(a_2)\}\\
&&-\frac{1}{\lambda} \{\lambda\varepsilon(a)\a(b_1)\o \a(b_2)+\nu\varepsilon(b)\a(a_1)\o \a(a_2)-\lambda\a^{2}(a)\o \a^{2}(b)\}\\
&=&\varepsilon(a)\a(b_1)\o \a(b_2)+\frac{\nu}{\lambda}\varepsilon(b)\a(a_1)\o \a(a_2)-\varepsilon(a)\a(b_1)\o \a(b_2)\\
&&-\frac{\nu}{\lambda}\varepsilon(b)\a(a_1)\o \a(a_2)+\a^{2}(a)\o \a^{2}(b)\\
&=&a\o b.
\end{eqnarray*}
It follows that $B\circ B^{-1}=id_{C\o C}$.
 Similarly, one may check that $B^{-1}\circ B=id_{C\o C}$.  So  $B^{-1}$ is the inverse of $B$.
$\hfill \Box$
\medskip

{\bf Example 3.3.}
Let $\{1,a,a^{2}\}$ be a basis of a 3-dimensional linear space $C$. The
following comultiplication $\Delta$,  the counit  $\varepsilon$ and the twist map $\alpha$ on $C$ define a Hom-coalgebra (\cite{yan}):
\begin{eqnarray*}
&&\Delta(1)=1\o 1, ~\Delta(a)=a^{2}\o a^{2},~ \Delta(a^{2})=a\o a,\\
&&\varepsilon(1)=1,~\varepsilon(a)=1,~\varepsilon(a^{2})=1,\\
&&\a(1)=1,~\a(a)=a^{2},~\a(a)=a^{2}.
\end{eqnarray*}
Therefore, by Theorem 3.1,  the solution $B$ of the HYBE for the Hom-coalgebra $C$ satisfies:
 \begin{eqnarray*}
&&B(1\o 1)=\nu 1\o 1,
~B(1\o a)=\lambda a^{2}\o a^{2}+\nu 1\o 1-\lambda 1\o a^{2},\\
&&B(1\o a^{2})=\lambda a\o a+\nu 1\o 1-\lambda 1\o a,
~B(a\o 1)=\lambda 1\o 1+\nu a^{2}\o a^{2}-\lambda a^{2}\o 1,\\
&&B(a\o a)=\nu a^{2}\o a^{2},
~B(a\o a^{2})=\lambda a\o a+ \nu a^{2}\o a^{2}-\lambda a^{2}\o a,\\
&&B(a^{2}\o 1)=\lambda 1\o 1+\nu a\o a-\lambda a\o 1,
~B(a^{2}\o a)=\nu a\o a,\\
&&B(a^{2}\o a^{2})=\lambda a^{2}\o a^{2}+\nu a\o a-\lambda a\o a^{2},
  \end{eqnarray*}
where $\lambda,\nu$ are two parameters in $k$.
 \medskip

\noindent{\bf Theorem 3.4.}
Let   $(C,\Delta, \varepsilon, \alpha)$ be a Hom-coalgebra and $\lambda,\nu\in k$.
Then
\begin{eqnarray*}
B: C\o C\rightarrow C\o C,~a\o b\mapsto \lambda \varepsilon(a)b_1\o b_2+\nu \varepsilon(b)a_1\o a_2-\nu \a(a)\o \a(b)
\end{eqnarray*}
is also a solution for HYBE.
If, in addition, $\lambda, \nu\in k^{\ast}$ and $\a$ is involutive, then  $B$  is   invertible,
where the inverse is given by
\begin{eqnarray*}
B^{-1}:C\o C\rightarrow C\o C,~a\o b\mapsto \frac{1}{\nu }\varepsilon(a)b_1\o b_2+\frac{1}{\lambda}\varepsilon(b)a_1\o a_2-\frac{1}{\nu} \a(a)\o \a(b).
\end{eqnarray*}

{\bf Proof.} Similar to the proof of Theorem 3.1 and Corollary 3.2.$\hfill \Box$
\medskip

\noindent{\bf Example 3.5.}
Let $\{1,g,x,y\}$ be a basis of a 4-dimensional linear space $H_{4}$. The
following comultiplication $\Delta$, counit $\varepsilon$ and the twist map $\alpha$ on $H_{4}$ define a Hom-coalgebra (\cite{Lu}):
\begin{eqnarray*}
&&\Delta(1)=1\o 1, \Delta(g)=g\o g,  \Delta(x)=kx\o 1+g\o kx, \Delta(y)=ky\o g+1\o ky, \\
&&\varepsilon(1)=1,~\varepsilon(g)=1,\varepsilon(x)=0,\varepsilon(y)=0,\\
&&\alpha(1)=1,~\alpha(g)=g,~\alpha(x)=kx,~\alpha(y)=ky,
\end{eqnarray*}
where $k$ is a parameter  in $k$.
By Theorem 3.4,  the solution $B$ of the HYBE for the Hom-algebra $H_{4}$ satisfies:
 \begin{eqnarray*}
&&B(1\o 1)=\lambda 1\o 1,
~B(1\o g)=\lambda g\o g+\nu 1\o 1-\nu 1\o g,\\
&&B(1\o x)=\lambda k x\o 1+\lambda kg\o x-\nu k 1\o x,
~B(1\o y)=\lambda k y\o g+\lambda k 1\o y-\nu k 1\o y,\\
&&B(g\o 1)=\lambda 1\o 1+\nu g\o g-\nu g\o 1,
~B(g\o g)=\lambda g\o g,\\
&&B(g\o x)=\lambda k x\o 1+\lambda k g\o x-\nu k g\o x,
~B(g\o y)=\lambda k y\o g+\nu k 1\o y-\nu k g\o y,\\
&&B(x\o 1)=\nu k g\o x,
~B(x\o g)=\nu k x\o 1+\nu k g\o x-\nu k x\o g,\\
&&B(x\o x)=-k^{2} x\o x,
~B(x\o y)=-k^{2} x\o y,
~B(y\o 1)=\nu k y\o g+\nu k 1\o y-\nu k y\o 1,\\
&&B(y\o g)=\nu k y\o g+\nu k 1\o y-\nu k y\o g,
~B(y\o x)=-k^{2} y\o x,
~B(y\o y)=-k^{2} y\o y,
  \end{eqnarray*}
where $\lambda,\nu$ are two parameters in $k$.

\section{Solutions of the HYBE from Hom-Lie algebras}
\def\theequation{\arabic{section}.\arabic{equation}}
\setcounter{equation} {0}

In this section, we obtain  a new solution of the HYBE from Hom-Lie algebras and prove that this solution is self-inverse.
\medskip

\noindent{\bf Theorem 4.1.}
Let  $(L, [\c, \c], \alpha)$  be a Hom-Lie algebra, $u$ is an $\a$-invariant element in $Z(L)$ and $\lambda, \upsilon\in k$.
Then
\begin{eqnarray*}
B: L\o L\rightarrow L\o L,~x\o y\mapsto \lambda [x,y]\o u- \upsilon\a(y)\o \a(x)
\end{eqnarray*}
is a solution for HYBE, where $Z(L)=\{u\in L|[u,x]=0, \forall x\in L\}$.
\medskip

{\bf Proof.} Obviously, $B$ is compatible with the twist map $\a$ since $z$ is   $\a$-invariant.
 Now we  verify that $B$ satisfy the HYBE. In fact, for any  $x,y,z\in L$,
 one may directly check that
$
 B(x\o u)=-\upsilon u\o \a(x), ~B(u\o x)=-\nu\a(x)\o u.
$
 On the one hand, we have
\begin{small}
\begin{eqnarray*}
&&(\a\o B)(x\o y\o z)\\
&&~~~~=\lambda\a(x)\o [y,z]\o u-\nu\a(x)\o\a(z)\o \a(y),\\
&&(B\o \a)\circ(\a\o B)(x\o y\o z)\\
&&~~~~=\lambda B(\a(x)\o [y,z])\o u-\nu B(\a(x)\o\a(z))\o \a^{2}(y)\\
&&~~~~=\lambda^{2}[\a(x), [y,z]]\o u\o u-\lambda\nu[\a(y),\a(z)] \o\a^{2}(x)\o u\\
&&~~~~~~~-\lambda\nu[\a(x),\a(z)] \o u\o \a^{2}(y)+\nu^2\a^{2}(z)\o \a^{2}(x)\o \a^{2}(y),\\
&&(\a\o B)\circ(B\o \a)\circ(\a\o B)(x\o y\o z)\\
&&~~~~=\lambda^{2}[\a^{2}(x), [\a(y),\a(z)]]\o B(u\o u)-\lambda\nu[\a^{2}(y),\a^{2}(z)] \o B(\a^{2}(x)\o u)\\
&&~~~~~~~-\lambda\nu[\a^{2}(x),\a^{2}(z)] \o B( u\o \a^{2}(y))+\nu^2\a^{3}(z)\o B(\a^{2}(x)\o \a^{2}(y))\\
&&~~~~=\underline{-\lambda^{2}\nu[\a^{2}(x), [\a(y),\a(z)]]\o u\o u}+\lambda\nu^{2}[\a^{2}(y),\a^{2}(z)] \o u\o\a^{3}(x) \\
&&~~~~~~~+\lambda\nu^{2}[\a^{2}(x),\a^{2}(z)] \o\a^{3}(y)\o u+\lambda\nu^{2}\a^{3}(z)\o[\a^{2}(x), \a^{2}(y)]\o u\\
&&~~~~~~~-\nu^{3}\a^{3}(z)\o\a^{3}(y)\o\a^{3}(x).
\end{eqnarray*}
\end{small}
On the other hand, we have
\begin{small}
\begin{eqnarray*}
&&(B\o \a)(x\o y\o z)\\
&&~~~~=\lambda [x,y]\o u\o \a(z)-\nu \a(y)\o \a(x)\o \a(z),\\
&&(\a\o B)\circ(B\o \a)(x\o y\o z)\\
&&~~~~=\lambda [\a(x),\a(y)]\o B(u\o \a(z))- \nu\a^{2}(y)\o B(\a(x)\o \a(z))\\
&&~~~~=-\lambda\nu[\a(x),\a(y)]\o\a^{2}(z)\o u-\lambda\nu\a^{2}(y)\o[\a(x), \a(z)]\o u\\
&&~~~~~~~+\nu^{2}\a^{2}(y)\o\a^{2}(z)\o\a^{2}(x),\\
&&(B\o \a)\circ(\a\o B)\circ(B\o \a)(x\o y\o z)\\
&&~~~~=-\lambda\nu B([\a(x),\a(y)]\o\a^{2}(z))\o u-\lambda\nu B(\a^{2}(y)\o[\a(x), \a(z)])\o u\\
&&~~~~~~~+\nu^{2}B(\a^{2}(y)\o\a^{2}(z))\o\a^{3}(x)\\
&&~~~~=\underline{-\lambda^{2}\nu[[\a(x),\a(y)],\a^{2}(z)] \o u\o u}+\lambda\nu^{2}\a^{3}(z)\o[\a^{2}(x),\a^{2}(y)]\o u\\
&&~~~~~~~\underline{-\lambda^{2}\nu[\a^{2}(y),[\a(x), \a(z)]] \o u\o u}+\lambda\nu^{2}[\a^{2}(x),\a^{2}(z)]\o\a^{3}(y)\o u\\
&&~~~~~~~+\lambda\nu^{2}[\a^{2}(y),\a^{2}(z)]\o u\o\a^{3}(x)-\nu^{3}\a^{3}(z)\o\a^{3}(y)\o\a^{3}(x).
\end{eqnarray*}
\end{small}
According to the anti-symmetry and the Hom-Jacobi identity, we have
\begin{eqnarray*}
[\a^{2}(x), [\a(y),\a(z)]]=[[\a(x),\a(y)],\a^{2}(z)]+[\a^{2}(y),[\a(x), \a(z)]].
\end{eqnarray*}
It follows that
\begin{eqnarray*}
(\a\o B)\circ(B\o \a)\circ(\a\o B)(x\o y\o z)=(B\o \a)\circ(\a\o B)\circ(B\o \a)(x\o y\o z).
\end{eqnarray*}
That is, $B$ is a solution for HYBE. $\hfill \Box$
\medskip

{\bf Corollary 4.2.}
Let  $(L, [\c, \c], \alpha)$  be a Hom-Lie algebra, $u$ is an $\a$-invariant element in $Z(L)$ and $\lambda\in k$.
Then
\begin{eqnarray*}
B: L\o L\rightarrow L\o L,~x\o y\mapsto \lambda [x,y]\o u- \a(y)\o \a(x)
\end{eqnarray*}
is a solution for HYBE.
If, in addition,  $\a$ is involutive, then $B$ is also invertible,
where the inverse is given by
\begin{eqnarray*}
B^{-1}: L\o L\rightarrow L\o L,~x\o y\mapsto \lambda u\o[x,y]- \a(y)\o \a(x).
\end{eqnarray*}
Furthermore, $B^{-1}$ is also a solution for HYBE.
\medskip

{\bf Proof.}
First, it is easy to see that $B$ is a solution for HYBE by setting $\nu=1$.

Next,  show that $B\circ B^{-1}=id_{L\o L}$. In fact, for any  $x,y,z\in L$,
we have
\begin{eqnarray*}
B\circ B^{-1}(x\o y)
&=&\lambda B(u\o[x,y])- B(\a(y)\o \a(x))\\
&=&-\lambda[\a(x), \a(y)]\o u-\lambda[\a(y), \a(x)]\o u+\a^{2}(y)\o \a^{2}(x)\\
&=&x\o y.
\end{eqnarray*}
It follows that $B\circ B^{-1}=id_{L\o L}$.
 Similarly, one may check that $B^{-1}\circ B=id_{L\o L}$.  So  $B^{-1}$ is the inverse of $B$.
Further,  similar to the proof of Theorem 4.1, one may calculate
\begin{eqnarray*}
&&(\a\o B^{-1})\circ(B^{-1}\o \a)\circ(\a\o B^{-1})(x\o y\o z)\\
&=&-\lambda^{2} u\o u \o[[\a(x),\a(y)],\a^{2}(z)]+\lambda\a^{3}(z)\o u\o[\a^{2}(x),\a^{2}(y)]\\
&&+\lambda u\o\a^{3}(y)\o[\a^{2}(x),\a^{2}(z)]+\lambda u\o[\a^{2}(y),\a^{2}(z)]\o\a^{3}(x)\\
&&-\a^{3}(z)\o\a^{3}(y)\o\a^{3}(x)\\
&=&(B^{-1}\o \a)\circ(\a\o B^{-1})\circ(B^{-1}\o \a)(x\o y\o z).
\end{eqnarray*}
So $B^{-1}$ is a solution for HYBE.
\medskip

{\bf Example 4.3.}
Let $(L,[\c,\c],\a)$ be a Hom-Lie algebra on 3-dimensional Euclidean $E^3$   with basis elements $\{e_{1},e_{2},e_{3}\}$,
 whose bracket $[\c,\c]$ is given by
 \begin{eqnarray*}
[e_{1},e_{2}]=e_{1},~ [e_{1},e_{3}]=0, ~[e_{2},e_{3}]=0.
  \end{eqnarray*}
The twist map $\alpha$ is given by
 \begin{eqnarray*}
 \alpha(e_{1})=e_{1}, ~\alpha(e_{2})=e_{2},~\alpha(e_{3})=-e_{3}.
  \end{eqnarray*}
Obviously, $e_{3}\in Z(L)$.
Therefore, by Theorem 4.1,  the solution $B$ of the HYBE for the Hom-Lie algebra $L$ satisfies:
 \begin{eqnarray*}
&&B(e_{1}\o e_{1})=-\nu e_{1}\o e_{1},
~B(e_{1}\o e_{2})=\lambda e_{1}\o e_{3}-\nu e_{2}\o e_{1},
~B(e_{1}\o e_{3})=\nu e_{3}\o e_{1},\\
&&B(e_{2}\o e_{1})=-\lambda e_{1}\o e_{3}-\nu e_{1}\o e_{2},
~B(e_{2}\o e_{2})=-\nu e_{2}\o e_{2},
~B(e_{2}\o e_{3})=\nu e_{3}\o e_{2},\\
&&B(e_{3}\o e_{1})=\nu e_{1}\o e_{3},
~B(e_{3}\o e_{2})=\nu e_{2}\o e_{3},
~B(e_{3}\o e_{3})=-\nu e_{3}\o e_{3},
  \end{eqnarray*}
where $\lambda,\nu \in k$.
\medskip

In the final part of this section, we present a kind of solutions of CHYBE from Hom-Lie algebras.
Recall from \cite{Yau2015},  Yau defined the CHYBE in a Hom-Lie algebra $(L, [\c,\c],\a)$ as
\begin{eqnarray*}
[r^{12},r^{13}]+[r^{12},r^{23}]+[r^{13},r^{33}]=0,
\end{eqnarray*}
for $r\in L^{\o 2}$. Here the three brackets above are defined as
 \begin{eqnarray*}
&&[r^{12},r^{13}]=[a_i, a_j]\o\a(b_i)\o\a(b_j),\\
&&[r^{12},r^{23}]=\a(a_i)\o[b_i, a_k]\o\a(a_k),\\
&&[r^{13},r^{23}]=\a(a_j)\o\a(a_k)\o[b_j, j_k],
\end{eqnarray*}
where $r^{12}=r\o 1=a_i\o b_i\o 1,r^{13}=(\tau\o id)(1\o r)=a_j\o 1\o b_j, $
$r^{23}=1\o r=1\o a_k\o b_k$.
\medskip

\noindent{\bf Theorem 4.4.}
Let  $(L, [\c, \c], \alpha)$  be a Hom-Lie algebra and  $u$ be an  element in $Z(L)$.
Then for any $x,y\in L $ and $m,n\in Z$,
\begin{eqnarray*}
r=\a^{m}([x,y])\o \a^{n}(u)
\end{eqnarray*}
is a solution for CHYBE.
\medskip

{\bf Proof.} It is easy to see that $[r^{12},r^{13}]=[r^{12},r^{23}]=[r^{13},r^{23}]=0$ since $u$ is an element in $Z(L)$.

\section{Hom-Yang-Baxter systems}
\def\theequation{\arabic{section}.\arabic{equation}}
\setcounter{equation} {0}

In this section, we extend the notion of Yang-Baxter systems to Hom-Yang-Baxter systems and present two kinds of Hom-Yang-Baxter systems.
\medskip

Consider three vector spaces $V, V', V''$,
let $\a_V, \a_{V'}, \a_{V''}$  be three endomorphisms on $V,V',V''$ and
$R: V\o V'\rightarrow V\o V', S: V\o V''\rightarrow V\o V'', T: V'\o V''\rightarrow V'\o V''$ be three linear maps.
Then a Hom-Yang-Baxter commutator is a map $[R, S, T]: V\o V'\o V'' \rightarrow  V\o V'\o V''$ defined by
\begin{eqnarray*}
[R, S, T]=R^{12}\circ S^{13}\circ T^{23}-T^{23}\circ S^{13}\circ R^{12},
\end{eqnarray*}
where $R^{12}=R\o \a_{V''},~S^{13}=(\tau\o id)\circ(\a_{V'}\o T)\circ (\tau\o id), ~T^{23}= \a_V\o T$.
\medskip

\noindent{\bf Definition 5.1.}
 Let $V, V'$ be two vector spaces and $\a_V, \a_{V'}$  be two endomorphisms.
 A system of linear maps
 \begin{eqnarray*}
W: V\o V\rightarrow V\o V, ~Z: V'\o V'\rightarrow V'\o V', ~X: V\o V'\rightarrow V\o V'
\end{eqnarray*}
is called a Hom-Yang-Baxter system, if the following conditions are satisfied:
 \begin{eqnarray}
[W, W, W]=0,~[Z, Z, Z]=0,~[W, X, X]=0,~[X, X, Z]=0.
\end{eqnarray}

\noindent{\bf Theorem 5.2.}
Let  $(A,\mu, 1_A, \alpha)$  be a Hom-algebra and $\lambda,\nu\in k$.
Then the following is a Hom-Yang-Baxter system:
\begin{eqnarray*}
&&W: A\o A\rightarrow A\o A,~a\o b\mapsto  ab\o 1_A+\lambda 1_A\o ab-\a(b)\o \a(a),\\
&&Z: A\o A\rightarrow A\o A,~a\o b\mapsto  \nu ab\o 1_A+1_A\o ab-\a(b)\o \a(a),\\
&&X: A\o A\rightarrow A\o A,~a\o b\mapsto  ab\o 1_A+1_A\o ab-\a(b)\o \a(a).
\end{eqnarray*}

{\bf Proof.}
It is sufficient to prove that the four equalities in Eq. (5.1) hold.
Here we only verify the equality $[W, X, X]=0$ and similar for other three equalities.
In fact, for any $a,b,c\in A$, on the one side, we have
\begin{small}
\begin{eqnarray*}
&&W^{12}\circ X^{13}\circ X^{23}(a\o b\o c)\\
&&~~~~=W^{12}\circ X^{13}\{\a(a)\o bc\o 1_A+\a(a)\o 1_A\o bc-\a(a)\o\a(b)\o\a(c)\}\\
&&~~~~=W^{12}\{\a^{2}(a)\o \a(bc)\o 1_A+\a(a)(bc)\o 1_A\o 1_A+1_A\o 1_A\o\a(a)(bc)\\
&&~~~~~~~-\a(bc)\o 1_A\o\a^{2}(a)-\a(ab)\o\a^{2}(c)\o 1_A-1_A\o\a^{2}(c)\o \a(ab)\\
&&~~~~~~~+a^{2}(b)\o a^{2}(c)\o a^{2}(a)\}\\
&&~~~~=\underline{\a^{2}(a)\a(bc)\o 1_A\o 1_A}_{(1)}
       +\underline{\lambda 1_A\o\a^{2}(a)\a(bc)\o 1_A}_{(2)}
       -\underline{\a^{2}(bc)\o \a^{3}(a)\o 1_A}_{(3)}\\
&&~~~~~~~+\underline{\a^{2}(a)\a(bc)\o 1_A\o 1_A}_{(4)}
         +\underline{\lambda 1_A\o\a^{2}(a)\a(bc)\o 1_A}_{(5)}
         - \underline{1_A\o \a^{2}(a)\a(bc)\o1_A}_{(6)}\\
&&~~~~~~~+\underline{\lambda  1_A\o 1_A\o\a^{2}(a)\a(bc)}_{(7)}
           -\underline{\a^{2}(bc)\o 1_A\o \a^{3}(a)}_{(8)}
           -\underline{\lambda 1_A\o \a^{2}(bc)\o \a^{3}(a)}_{(9)}\\
&&~~~~~~~+\underline{1_A\o \a^{2}(bc)\o \a^{3}(a)}_{(10)}
         -\underline{\a(ab)\a^{2}(c)\o 1_A\o 1_A}_{(11)}
         -\underline{\lambda 1_A\o \a(ab)\a^{2}(c)\o 1_A}_{(12)}\\
&&~~~~~~~+\underline{\a^{3}(c)\o \a^{2}(ab)\o 1_A}_{(13)}
         -\underline{\lambda 1_A\o\a^{3}(c)\o \a^{2}(ab)}_{(14)}
         +\underline{\a^{2}(bc)\o1_A\o\a^{3}(a)}_{(15)}\\
&&~~~~~~~+\underline{\lambda  1_A\o  \a^{2}(bc)\o \a^{3}(a)} _{(16)}
         -\underline{\a^{3}(c)\o\a^{3}(b)\o\a^{3}(a)}_{(17)}\\
&&~~~~=\underline{\a^{2}(a)\a(bc)\o 1_A\o 1_A}_{(1)}
       +\underline{\lambda 1_A\o\a^{2}(a)\a(bc)\o 1_A}_{(2)}
       -\underline{\a^{2}(bc)\o \a^{3}(a)\o 1_A}_{(3)}\\
&&~~~~~~~-\underline{1_A\o \a^{2}(a)\a(bc)\o1_A}_{(6)}
          +\underline{\lambda  1_A\o 1_A\o\a^{2}(a)\a(bc)}_{(7)}
          +\underline{1_A\o \a^{2}(bc)\o \a^{3}(a)}_{(10)}\\
&&~~~~~~~+\underline{\a^{3}(c)\o \a^{2}(ab)\o 1_A}_{(13)}
         -\underline{\lambda 1_A\o\a^{3}(c)\o \a^{2}(ab)}_{(14)}
         -\underline{\a^{3}(c)\o\a^{3}(b)\o\a^{3}(a)}_{(17)}.
\end{eqnarray*}
\end{small}
The last  equality holds since (4)=(11), (5)=(12), (8)=(15) and  (9)=(16).

 On  the other side, we have
 \begin{small}
\begin{eqnarray*}
&&X^{23}\circ X^{13}\circ W^{12}(a\o b\o c)\\
&&~~~~=X^{23}\circ X^{13}\{ab\o 1_A\o \a(c)+\lambda 1_A\o ab\o \a(c)-\a(b)\o \a(a)\o \a(c)\}\\
&&~~~~=X^{23}\{(ab)\a(c)\o 1_A\o 1_A+1_A\o 1_A\o(ab)\a(c)-\a^{2}(c)\o 1_A\o\a(ab)\\
&&~~~~~~~+\lambda 1_A\o\a(ab)\o\a^{2}(c)-\a(bc)\o\a^{2}(a)\o 1_A-1_\o\a^{2}(a)\o\a(bc)\\
&&~~~~~~~+\a^{2}(c)\o\a^{2}(a)\o\a^{2}(b)\}\\
&&~~~~=\underline{\a(ab)\a^{2}(c)\o 1_A\o 1_A}_{(18)}
      +\underline{1_A\o 1_A\o\a(ab)\a^{2}(c)}_{(19)}
      -\underline{\a^{3}(c)\o 1_A\o\a^{2}(ab)}_{(20)}\\
&&~~~~~~~+\underline{\lambda 1_A\o\a(ab)\o\a^{2}(c)}_{(21)}
         +\underline{\lambda 1_A\o 1_A\o\a(ab)\a^{2}(c)}_{(22)}
         -\underline{\lambda 1_A\o\a^{3}(c)\o\a^{2}(ab)}_{(23)}\\
&&~~~~~~~-\underline{a^{2}(bc)\o \a^{3}(a)\o 1_A }_{(24)}
         -\underline{1_A\o\a^{2}(a)\a(bc)\o 1_A }_{(25)}
         -\underline{1_A\o1_A\o\a^{2}(a)\a(bc)}_{(26)}\\
&&~~~~~~~+\underline{1_A\o \a^{2}(bc)\o  \a^{3}(a)}_{(27)}
          +\underline{\a^{3}(c)\o\a^{2}(ab)\o1_A}_{(28)}
          +\underline{\a^{3}(c)\o1_A  \o\a^{2}(ab)}_{(29)} \\
&&~~~~~~~ -\underline{\a^{3}(c)\o\a^{3}(b)\o\a^{3}(a)}_{(30)} \\
&&~~~~=\underline{\a(ab)\a^{2}(c)\o 1_A\o 1_A}_{(18)}
        +\underline{\lambda 1_A\o\a(ab)\o\a^{2}(c)}_{(21)}
         +\underline{\lambda 1_A\o 1_A\o\a(ab)\a^{2}(c)}_{(22)}\\
&&~~~~~~~-\underline{\lambda 1_A\o\a^{3}(c)\o\a^{2}(ab)}_{(23)}
          -\underline{a^{2}(bc)\o \a^{3}(a)\o 1_A }_{(24)}
         -\underline{1_A\o\a^{2}(a)\a(bc)\o 1_A }_{(25)}\\
 &&~~~~~~~+\underline{1_A\o \a^{2}(bc)\o  \a^{3}(a)}_{(27)}
          +\underline{\a^{3}(c)\o\a^{2}(ab)\o1_A}_{(28)}
           -\underline{\a^{3}(c)\o\a^{3}(b)\o\a^{3}(a)}_{(30)}.
\end{eqnarray*}
\end{small}
The last  equality holds since (19)=(26) and  (20)=(29).

Compare the two expressions above, we have
 \begin{eqnarray*}
 &&(1)=(18), (2)=(21), (3)=(24), (6)=(25), \\
 &&(7)=(22), (10)=(27), (13)=(28), (14)=(23), (17)=(30).
 \end{eqnarray*}
It follows that $[W, X, X]=0$, as desired.
The proof is finished. $\hfill \Box$
\medskip

\noindent{\bf Theorem 5.3.}
Let  $(C,\Delta, \varepsilon, \alpha)$  be a Hom-coalgebra and $\lambda,\nu\in k$.
Then the following is a Hom-Yang-Baxter system:
\begin{eqnarray*}
&&W: C\o C\rightarrow C\o C,~a\o b\mapsto \lambda \varepsilon(a)b_1\o b_2+\varepsilon(b)a_1\o a_2- \a(b)\o \a(a),\\
&&Z: C\o C\rightarrow C\o C,~a\o b\mapsto \varepsilon(a)b_1\o b_2+\nu \varepsilon(b)a_1\o a_2-\a(b)\o \a(a),\\
&&X: C\o C\rightarrow C\o C,~a\o b\mapsto  \varepsilon(a)b_1\o b_2+\varepsilon(b)a_1\o a_2-\a(b)\o \a(a).
\end{eqnarray*}

{\bf Proof.}
We only prove $[X, X, Z]=0$ and similar for other three equalities.
In fact, for any $a,b,c\in C$, on the one side, we have
\begin{small}
\begin{eqnarray*}
&&Z^{23}\circ X^{13}\circ X^{12}(a\o b\o c)\\
&&~~~~=Z^{23}\circ X^{13}\{\varepsilon(a)b_1\o b_2\o\a(c)+\nu \varepsilon(b)a_1\o a_2\o\a(c)-\a(b)\o \a(a)\o\a(c)\}\\
&&~~~~=Z^{23}\{\varepsilon(a)\a(c_1)\o \a^{2}(b)\o\a(c_2)+\varepsilon(a)\varepsilon(c)b_{11}\o\a(b_2)\o b_{12}\\
&&~~~~~~~-\varepsilon(a)\a^{2}(c)\o\a(b_2)\o\a(b_1)+\underline{\varepsilon(b)\a(c_1)\o \a^{2}(a)\o\a(c_2)}\\
&&~~~~~~~+\varepsilon(b)\varepsilon(c)a_{11}\o\a(a_2)\o a_{12}-\varepsilon(b)\a^{2}(c)\o\a(a_2)\o\a(a_1)\\
&&~~~~~~~-\underline{\varepsilon(b)\a(c_1)\o \a^{2}(a)\o\a(c_2)}-\varepsilon(c)\a(b_1)\o \a^{2}(a)\o\a(b_2)\\
&&~~~~~~~+\a^{2}(c)\o\a^{2}(a)\o\a^{2}(b)\}\\
&&~~~~=\underline{\varepsilon(a)\varepsilon(b)\a^{2}(c_{1})\o \a(c_{21})\o \a(c_{22})}_{(1)}
       +\underline{\nu\varepsilon(a)\a^{3}(c)\o\a^{2}(b_{1})\o \a^{2}(b_{2})}_{(2)}\\
&&~~~~~~~-\underline{\varepsilon(a)\a^{2}(c_{1})\o \a^{2}(c_{2})\o\a^{3}(b)}_{(3)}
       + \underline{\varepsilon(a)\varepsilon(c)\a^{2}(b_{1})\o \a(b_{21})\o \a(b_{22})}_{(4)}\\
&&~~~~~~~+\underline{\nu\varepsilon(a)\varepsilon(c)\a^{2}(b_{1})\o \a(b_{21})\o \a(b_{22})}_{(5)}
         - \underline{\varepsilon(a)\varepsilon(c)\a(b_{11})\o \a(b_{12}) \o\a^{2}(b_{2})}_{(6)}\\
&&~~~~~~~-\underline{\varepsilon(a)\a^{3}(c)\o\a^{2}(b_{1})\o \a^{2}(b_{2})}_{(7)}
         -\underline{\nu\varepsilon(a)\a^{3}(c)\o\a^{2}(b_{1})\o \a^{2}(b_{2})}_{(8)}\\
&&~~~~~~~+\underline{\varepsilon(a)\a^{3}(c)\o\a^{2}(b_{1})\o \a^{2}(b_{2})  }_{(9)}
         +\underline{\varepsilon(b)\varepsilon(c)\a^{2}(a_{1})\o \a(a_{21})\o \a(a_{22})}_{(10)} \\
&&~~~~~~~+\underline{\nu\varepsilon(b)\varepsilon(c)\a^{2}(a_{1})\o \a(a_{21})\o \a(a_{22})}_{(11)}
         -\underline{\varepsilon(b)\varepsilon(c)\a(a_{11})\o \a(a_{12}) \o\a^{2}(a_{2})}_{(12)}\\
&&~~~~~~~- \underline{\varepsilon(b) \a^{3}(c)\o\a^{2}(a_{1})\o \a^{2}(a_{2}) }_{(13)}
          -\underline{\nu \varepsilon(b) \a^{3}(c)\o\a^{2}(a_{1})\o \a^{2}(a_{2})}_{(14)}\\
&&~~~~~~~+\underline{\varepsilon(b) \a^{3}(c)\o\a^{2}(a_{1})\o \a^{2}(a_{2})  }_{(15)}
         -\underline{\varepsilon(a)\varepsilon(c)\a^{2}(b_{1})\o \a(b_{21})\o \a(b_{22})}_{(16)}\\
&&~~~~~~~-\underline{\nu \varepsilon(c) \a^{3}(b)\o\a^{2}(a_{1})\o \a^{2}(a_{2}) }_{(17)}
         +\underline{\varepsilon(c)\a^{2}(b_{1})\o \a^{2}(b_{2})\o \a^{3}(a)}_{(18)} \\
&&~~~~~~~+\underline{\varepsilon(a)\a^{3}(c)\o\a^{2}(b_{1})\o \a^{2}(b_{2}) }_{(19)}
          +\underline{\nu\varepsilon(b)\a^{3}(c)\o\a^{2}(a_{1})\o \a^{2}(a_{2})}_{(20)}  \\
&&~~~~~~~-  \underline{\a^{3}(c)\o\a^{3}(b)\o\a^{3}(a)}_{(21)}\\
&&~~~~=\underline{\varepsilon(a)\varepsilon(b)\a^{2}(c_{1})\o \a(c_{21})\o \a(c_{22})}_{(1)}
       -\underline{\varepsilon(a)\a^{2}(c_{1})\o \a^{2}(c_{2})\o\a^{3}(b)}_{(3)}\\
&&~~~~~~~+\underline{\nu\varepsilon(a)\varepsilon(c)\a^{2}(b_{1})\o \a(b_{21})\o \a(b_{22})}_{(5)}
         +\underline{\varepsilon(a)\a^{3}(c)\o\a^{2}(b_{1})\o \a^{2}(b_{2})  }_{(9)} \\
&&~~~~~~~+\underline{\nu\varepsilon(b)\varepsilon(c)\a^{2}(a_{1})\o \a(a_{21})\o \a(a_{22})}_{(11)}
          -\underline{\varepsilon(a)\varepsilon(c)\a^{2}(b_{1})\o \a(b_{21})\o \a(b_{22})}_{(16)}\\
&&~~~~~~~-\underline{\nu \varepsilon(c) \a^{3}(b)\o\a^{2}(a_{1})\o \a^{2}(a_{2}) }_{(17)}
         +\underline{\varepsilon(c)\a^{2}(b_{1})\o \a^{2}(b_{2})\o \a^{3}(a)}_{(18)} \\
&&~~~~~~~-  \underline{\a^{3}(c)\o\a^{3}(b)\o\a^{3}(a)}_{(21)}.
\end{eqnarray*}
\end{small}
The last  equality holds since (2)=(8), (4)=(6), (7)=(9),(10)=(12),(13)=(15) and  (14)=(20).

On the other side, we have
\begin{small}
\begin{eqnarray*}
&&X^{12}\circ X^{13}\circ Z^{23}(a\o b\o c)\\
&&~~~~=X^{12}\circ X^{13}\{\varepsilon(b)\a(a)\o c_1\o c_2+\nu\varepsilon(c)\a(a)\o b_1\o b_2-\a(a)\o\a(c)\o\a(b)\}\\
&&~~~~=X^{12}\{\varepsilon(a)\varepsilon(b)c_{21}\o\a(c_1)\o c_{22}+\underline{\varepsilon(b)\a(a_1)\o\a^{2}(c)\o\a(a_2)}\\
&&~~~~~~~-\varepsilon(b)\a(c_2)\o\a(c_1)\o\a^{2}(a)+\nu\varepsilon(a)\varepsilon(c)b_{21}\o\a(b_1)\o b_{22}\\
&&~~~~~~~+\nu\varepsilon(c)\a(a_1)\o\a^{2}(b)\o\a(a_2)-\nu\varepsilon(c)\a(b_2)\o\a(b_1)\o\a^{2}(a)\\
&&~~~~~~~-\varepsilon(a)\a(b_1)\o\a^{2}(c)\o\a(b_2)-\underline{\varepsilon(b)\a(a_1)\o\a^{2}(c)\o\a(a_2)}\\
&&~~~~~~~+\a^{2}(b)\o\a^{2}(c)\o\a^{2}(a)\}\\
&&~~~~=\underline{\varepsilon(a)\varepsilon(b)\a(c_{11})\o\a(c_{12})\o \a^{2}(c_{2})}_{(22)}
       +\underline{\varepsilon(a)\varepsilon(b)\a(c_{11})\o\a(c_{12})\o \a^{2}(c_{2})}_{(23)}\\
&&~~~~~~~-\underline{\varepsilon(a)\varepsilon(b)\a(c_{11})\o\a(c_{12})\o \a^{2}(c_{2})}_{(24)}
        -\underline{\varepsilon(b) \a^{2}(c_{1})\o \a^{2}(c_{2})\o  \a^{3}(a)}_{(25)}\\
&&~~~~~~~- \underline{\varepsilon(b) \a^{2}(c_{1})\o \a^{2}(c_{2})\o  \a^{3}(a)}_{(26)}
        +\underline{\varepsilon(b) \a^{2}(c_{1})\o \a^{2}(c_{2})\o  \a^{3}(a)}_{(27)}\\
&&~~~~~~~+\underline{\nu\varepsilon(a)\varepsilon(c) \a(b_{11})\o\a(b_{12})\o \a^{2}(b_{2}) }_{(28)}
         +\underline{\nu\varepsilon(a)\varepsilon(c) \a(b_{11})\o\a(b_{12})\o \a^{2}(b_{2})}_{(29)}\\
&&~~~~~~~- \underline{\nu\varepsilon(a)\varepsilon(c) \a(b_{11})\o\a(b_{12})\o \a^{2}(b_{2})}_{(30)}
          +\underline{\nu\varepsilon(c) \a^{2}(b_{1})\o \a^{2}(b_{2})\o  \a^{3}(a)}_{(31)}\\
&&~~~~~~~+\underline{\nu\varepsilon(b)\varepsilon(c) \a(a_{11})\o\a(a_{12})\o \a^{2}(a_{2}) } _{(32)}
        -\underline{\nu\varepsilon(c) \a^{3}(b)\o\a^{2}(a_{1})\o \a^{2}(a_{2})}_{(33)}\\
&&~~~~~~~-\underline{\nu\varepsilon(c) \a^{2}(b_{1})\o \a^{2}(b_{2})\o  \a^{3}(a)}_{(34)}
        -\underline{\nu\varepsilon(c) \a^{2}(b_{1})\o \a^{2}(b_{2})\o  \a^{3}(a)}_{(35)}\\
&&~~~~~~~+\underline{\nu\varepsilon(c) \a^{2}(b_{1})\o \a^{2}(b_{2})\o  \a^{3}(a)}_{(36)}
         -\underline{\varepsilon(a) \a^{2}(c_{1})\o \a^{2}(c_{2})\o  \a^{3}(b)}_{(37)}\\
&&~~~~~~~- \underline{\varepsilon(a)\varepsilon(c)\a(b_{11})\o\a(b_{12})\o \a^{2}(b_{2}) } _{(38)}
         +\underline{\varepsilon(a)\a^{3}(c)\o\a^{2}(b_{1})\o \a^{2}(b_{2})}_{(39)}\\
&&~~~~~~~+\underline{\varepsilon(b)\a^{2}(c_{1})\o \a^{2}(c_{2})\o  \a^{3}(a)  } _{(40)}
          +\underline{\varepsilon(c) \a^{2}(b_{1})\o \a^{2}(b_{2})\o  \a^{3}(a)}_{(41)}\\
 &&~~~~~~~-  \underline{\a^{3}(c)\o\a^{3}(b)\o\a^{3}(a)}_{(42)}\\
&&~~~~= \underline{\varepsilon(a)\varepsilon(b)\a(c_{11})\o\a(c_{12})\o \a^{2}(c_{2})}_{(22)}
       +\underline{\nu\varepsilon(a)\varepsilon(c) \a(b_{11})\o\a(b_{12})\o \a^{2}(b_{2}) }_{(28)}\\
&&~~~~~~~+\underline{\nu\varepsilon(b)\varepsilon(c) \a(a_{11})\o\a(a_{12})\o \a^{2}(a_{2}) } _{(32)}
        -\underline{\nu\varepsilon(c) \a^{3}(b)\o\a^{2}(a_{1})\o \a^{2}(a_{2})}_{(33)}\\
&&~~~~~~~-\underline{\varepsilon(a) \a^{2}(c_{1})\o \a^{2}(c_{2})\o  \a^{3}(b)}_{(37)}
         - \underline{\varepsilon(a)\varepsilon(c)\a(b_{11})\o\a(b_{12})\o \a^{2}(b_{2}) } _{(38)}\\
&&~~~~~~~+\underline{\varepsilon(a)\a^{3}(c)\o\a^{2}(b_{1})\o \a^{2}(b_{2})}_{(39)}
          +\underline{\varepsilon(c) \a^{2}(b_{1})\o \a^{2}(b_{2})\o  \a^{3}(a)}_{(41)}\\
 &&~~~~~~~-  \underline{\a^{3}(c)\o\a^{3}(b)\o\a^{3}(a)}_{(42)}.
\end{eqnarray*}
\end{small}
The last  equality holds since (23)=(24), (25)=(40), (26)=(27),(29)=(30),(31)=(34) and  (35)=(36).
Compare the two expressions above, we have
 \begin{eqnarray*}
 &&(1)=(22), (3)=(37), (5)=(28), (9)=(39), \\
 &&(11)=(32), (17)=(33), (18)=(41), (16)=(38), (21)=(42).
 \end{eqnarray*}
It follows that $[X, X, Z]=0$, as desired.
The proof is finished. $\hfill \Box$
\medskip

\begin{center}
 {\bf ACKNOWLEDGEMENT}
 \end{center}
  The paper is  supported by the NSF of China (Nos. 11761017 and 11801304),
   Guizhou Provincial  Science and Technology  Foundation (No. [2020]1Y005),
    the Anhui Provincial Natural Science Foundation (No. 1908085MA03)
     and the Key University Science Research Project of Anhui Province  (No. KJ2020A0711).

{\bf Data Availability Statement} Our manuscripts does not include a data availability statement.

\renewcommand{\refname}{REFERENCES}

\end{document}